%% file: qci_restriction_persistent_nodal_curves.tex
\documentclass[11pt]{amsart}
\usepackage[a4paper,margin=26mm]{geometry}
\usepackage[T1]{fontenc}
\usepackage{lmodern,microtype}
\usepackage{amsmath,amssymb,mathtools}
\usepackage{graphicx,float,enumitem,booktabs}
\usepackage{xcolor}
\usepackage[colorlinks=true,linkcolor=blue!45!black,citecolor=blue!45!black,urlcolor=blue!45!black]{hyperref}
\numberwithin{equation}{section}
\newtheorem{theorem}{Theorem}[section]
\newtheorem{proposition}[theorem]{Proposition}
\newtheorem{corollary}[theorem]{Corollary}
\newtheorem{lemma}[theorem]{Lemma}
\theoremstyle{definition}

\theoremstyle{remark}
\newtheorem{remark}[theorem]{Remark}
\newcommand{\T}{\mathbb T}
\newcommand{\R}{\mathbb R}
\newcommand{\C}{\mathbb C}
\newcommand{\E}{\mathcal E}
\newcommand{\norm}[1]{\left\lVert #1\right\rVert}
\newcommand{\Op}{\operatorname{Op}^w_h}

\DeclareMathOperator{\supp}{supp}
\DeclareMathOperator{\Impart}{Im}
\DeclareMathOperator{\Repart}{Re}
\title[Restriction bounds and common nodal curves]{Uniform restriction lower bounds and nongeodesic common nodal curves in integrable systems}
\author{Xianchao Wu}
\address{School of Mathematics and Statistics, Wuhan University of Technology, Wuhan, Hubei, China}
\email{xianchao.wu@whut.edu.cn}
\author{Xiao Xiao}
\address{Department of Mathematics and Statistics, McGill University, Montr\'eal, Quebec, Canada}
\email{xiao.xiao5@mail.mcgill.ca}
\date{}
\begin{document}
\begin{abstract}
We construct a family of real analytic metrics on the two-torus, for which a fixed embedded nongeodesic circle is a nodal component of infinitely many real Laplace eigenfunctions. The examples address a question asked by Bourgain-Rudnick \cite{BR11}, and give a negative answer to the questions asked by Toth-Zelditch \cite{TZ21} and Han \cite{Han}: whether a curve shared by infinitely many nodal sets must be geodesic.

To explain the restriction obstruction, we prove uniform $L^2$ lower bounds for joint eigenfunctions of quantum completely integrable (QCI) systems under transverse crossing and isolation of a tangential characteristic branch, then extend to multi-branch setting. The results are used to characterize the possible joint nodal curves in the Liouville tori, narrowing down the choices. Finally, the counterexample realizes one possibility for which cancellation happen exactly on a nongeodesic curve.

\end{abstract}
\maketitle
\input{sections/introduction}
\input{sections/core}
\input{sections/applications}
\input{sections/finite_rigidity}
\input{sections/counterexample}

\input{sections/references}
\end{document}

%% file: sections/introduction.tex
\section{Introduction and main results}\label{sec:introduction}

A fixed curve can occur in the nodal sets of eigenfunctions at arbitrarily large eigenvalues. On a flat two-torus, a theorem of Bourgain--Rudnick says that a real analytic curve has this property precisely when it is a segment of a closed geodesic \cite[Theorem 1.1]{BR11}. Does the same geometric rigidity hold on a general analytic surface? Toth--Zelditch raise this question in their discussion of nodal curves \cite[\S1.3 of the arXiv version]{TZ21}; it is also formulated in Han's problem notes \cite[Problem 4.9]{Han}.

We give an explicit negative answer. One fixed analytic metric, with no nonzero Killing field, admits a fixed embedded nongeodesic circle which is a regular nodal component of an infinite sequence of real joint eigenfunctions. The construction persists within an explicit family converging to a flat metric. Thus analyticity or the absence of continuous symmetries do not guarantee that the common nodal curves is a geodesic segment.

\subsection{The nongeodesic counterexample}
Fix a metric and a curve. Call a fixed curve \emph{persistent} if it is contained in the nodal sets of eigenfunctions with eigenvalues tending to infinity. 

\begin{theorem}[A nongeodesic common nodal circle]\label{thm:counterexample}
On $M=(\R/2\pi\mathbb Z)^2$, let $0<\varepsilon<2$ and set
\begin{equation}\label{eq:explicit-metric}
 g_\varepsilon=F_\varepsilon(x,y)(dx^2+dy^2),\qquad
 F_\varepsilon=10+\varepsilon\cos3x+4\varepsilon\cos6y.
\end{equation}
There are distinct eigenvalues $\lambda_j\to\infty$ and real $L^2(M,g_\varepsilon)$-normalized eigenfunctions $u_j$ such that
\begin{equation}\label{eq:explicit-curve}
 \gamma(t)=(2t,t)\pmod{2\pi},\qquad t\in\R/2\pi\mathbb Z,
\end{equation}
is a regular connected component of every nodal set $\{u_j=0\}$. If
\[
 a(x)=2+\varepsilon\cos3x,\qquad b(y)=8+4\varepsilon\cos6y,
 \qquad Q_\varepsilon=F_\varepsilon^{-1}(-b\partial_x^2+a\partial_y^2),
\]
then $-\Delta_{g_\varepsilon}u_j=\lambda_j u_j$ and $Q_\varepsilon u_j=0$.
The metric has no nonzero Killing vector field, and the signed geodesic curvature of the oriented circle \eqref{eq:explicit-curve} is
\begin{equation}\label{eq:explicit-curvature}
 k_{g_\varepsilon}(t)=\frac{9\varepsilon\sin6t}{10(2+\varepsilon\cos6t)^{3/2}}.
\end{equation}
In particular, the common nodal circle is not a geodesic. Moreover, $g_\varepsilon\to10(dx^2+dy^2)$ in $C^\infty$ as $\varepsilon\downarrow0$.
\end{theorem}

The analytic part of the paper identifies the relevant obstruction. Restriction records the tangential component of a characteristic covector. If a branch with positive joint Liouville mass is isolated by this tangential projection and crosses the hypersurface transversely, a microlocal Cauchy estimate gives a uniform restriction lower bound. Distinct branches with the same tangential covector can instead cancel. Our example realizes this cancellation exactly on a nongeodesic curve. The general restriction theorem, its extension to disconnected fibers, and the explicit construction therefore give complementary information about persistent nodal sets.

The winding numbers of $\gamma$ are $(2,1)$. Formula \eqref{eq:explicit-curvature} distinguishes a coordinate straight line from a geodesic of the metric. Its curvature is nonzero on open arcs, so even such arcs can be shared by infinitely many nodal sets.

The proof in Section~\ref{sec:counterexample} uses nonreal Floquet multipliers for a one-dimensional weighted eigenvalue problem. The same scalar solutions, inserted at two different coordinate scales, give exact joint eigenfunctions with a common phase front. The restriction theorem below supplies necessary geometric conditions for persistence, while the construction shows that the nongeodesic alternative can occur.

\subsection{Restriction on a regular joint torus}
Let $M$ be a compact smooth manifold without boundary, of dimension $n$, equipped with a smooth positive density. Let
\[
 P(h)=(P_1(h),\ldots,P_n(h)),\qquad P_j(h)\in\Psi_h^1(M),
\]
be commuting self-adjoint operators whose full symbols have the form $p_j+h p_{j,1}(h)$, with $p_j$ real and the lower-order symbols uniformly bounded in the appropriate symbol classes. Assume that $\sum_jP_j(h)^2$ is uniformly semiclassically elliptic of order two. Write
\[
 p=(p_1,\ldots,p_n):T^*M\longrightarrow\R^n,
 \qquad \Lambda_E=p^{-1}(E).
\]
We assume that $E$ is a regular value, so $dp_1\wedge\cdots\wedge dp_n\ne0$ on the compact joint fiber. Each connected component is a Lagrangian torus by the Liouville--Arnold theorem. Consider normalized exact joint eigenfunctions
\begin{equation}\label{eq:general-joint}
 P_j(h)u_h=E_j(h)u_h,\qquad E(h)\longrightarrow E,
 \qquad \norm{u_h}_{L^2(M)}=1.
\end{equation}
The customary assumption $E(h)=E+O(h)$ is included. Allowing convergence without a specified rate is justified by the parameter-dependent local construction in Section~\ref{sec:cauchy}.

Let $H\subset M$ be a compact smooth embedded hypersurface, possibly with boundary; all patches below lie in its interior. For a local defining function $\rho$, set
\[
 T_H^*M=\{(x,\xi):x\in H\},\qquad
 \Sigma_E=\Lambda_E\cap T_H^*M,\qquad
 r_H(x,\xi)=(x,\xi|_{T_xH}).
\]
An \emph{isolated transverse point} $z_0\in\Sigma_E$ means the following. There is an index $j$ with $H_{p_j}\rho(z_0)\ne0$, a sufficiently small neighborhood $U$ of $z_0$, and an open set $W\subset T^*H$ containing $r_H(z_0)$ such that
\begin{equation}\label{eq:branch-isolation}
 \Sigma_E\cap r_H^{-1}(W)\subset U.
\end{equation}
Here $U$ is chosen so that, in coordinates $H=\{x_1=0\}$, $\partial_{\xi_1}p_j\ne0$ and $p_j=E_j$ is a single graph $\xi_1=a_0(x,\xi')$ in $U$. Thus \eqref{eq:branch-isolation} excludes competing points of the \emph{entire} joint fiber over the selected tangential window. This definition is independent of the choice of defining function and coordinates.

\begin{theorem}[Restriction on a connected regular joint fiber]\label{thm:connected}
Assume that $\Lambda_E$ is connected and that $H$ has an isolated transverse point on $\Lambda_E$. Then the joint eigenfunctions \eqref{eq:general-joint} satisfy
\begin{equation}\label{eq:main-lower}
 \norm{u_h|_H}_{L^2(H)}\ge c>0
\end{equation}
for all sufficiently small $h$.
\end{theorem}

In particular, the hypotheses hold if $r_H$ is injective on $\Sigma_E$ and some Hamiltonian flow crosses $H$ transversely there. The local formulation \eqref{eq:branch-isolation} is sufficient because a single open flow box has positive normalized Liouville measure on a connected joint torus. 

\subsection{Finite fibers and the geometry of persistence}
The rotating joint fiber of a Liouville torus has four connected components. To apply the restriction argument to this model, each component must have an observable patch.

\begin{theorem}[Restriction on finitely many joint tori]\label{thm:finite}
Suppose that $\Lambda_E=\bigsqcup_{\nu=1}^{N}\Lambda_{E,\nu}$ is a compact regular joint fiber. Assume that every component $\Lambda_{E,\nu}$ contains an isolated transverse point over $H$, with isolation always understood relative to the full fiber as in \eqref{eq:branch-isolation}. Then \eqref{eq:main-lower} holds for every normalized sequence \eqref{eq:general-joint}.
\end{theorem}

The components may carry unknown and sequence-dependent proportions of the mass. The proof in Section~\ref{sec:finite} uses the decomposition of a defect measure into normalized Liouville measures on the finitely many components.

For a Liouville metric $g=(a(x)+b(y))(dx^2+dy^2)$ and a rotating joint energy $\kappa$, write $A_\kappa=\sqrt{a+\kappa}$ and $B_\kappa=\sqrt{b-\kappa}$. The four tangential momenta at a nonzero curve tangent $v$ are
\[
 \sigma A_\kappa v_x+\tau B_\kappa v_y,\qquad \sigma,\tau\in\{+1,-1\}.
\]
They are distinct precisely when
\begin{equation}\label{eq:intro-separation}
 v_xv_y\ne0,\qquad A_\kappa^2v_x^2\ne B_\kappa^2v_y^2.
\end{equation}
Theorem~\ref{thm:finite} therefore gives a restriction lower bound near any such point. Its contrapositive forces every connected persistent joint nodal curve to follow one of the four line fields
\[
 \ker dx,\quad\ker dy,\quad
 \ker(A_\kappa dx+B_\kappa dy),\quad
 \ker(A_\kappa dx-B_\kappa dy).
\]
Section~\ref{sec:finite} proves this rigidity statement and its consequences: a fixed contractible closed curve cannot persist in a compact rotating window, a closed noncoordinate persistent curve requires an action resonance, and restriction on a curve with a separating point determines an exact joint eigenfunction stably.

The circle of Theorem~\ref{thm:counterexample} follows one of the permitted phase fronts. Its action resonance is $I_y(0)=2I_x(0)$. Along the circle, two distinct transverse characteristic branches have the same tangential covector and their traces cancel. This realizes the obstruction to Theorem~\ref{thm:finite}. The same local Cauchy estimate determines the quantum limit of the constructed real sequence: it is the equal mixture of Liouville measures on the two cancelling components (Proposition~\ref{prop:quantum-limit}). Because its full joint fiber has four components, we retain a separate spherical example in Section~\ref{sec:sphere} which demonstrates the necessity of branch isolation within the connected-fiber setting of Theorem~\ref{thm:connected}.

\subsection{Related results and organization}
Upper restriction bounds for Laplace eigenfunctions were established by Burq--G\'erard--Tzvetkov \cite{BGT07}; for a hypersurface the general bound is $O(h^{-1/4})$. In the quantum integrable surface setting, Toth obtained logarithmic bounds for suitable curves \cite{Tot09}. Bourgain--Rudnick proved uniform toral restriction estimates \cite{BR09,BR11}. Quantum ergodic restriction provides lower bounds along suitable subsequences under dynamical and asymmetry assumptions \cite{TZ13,DZ13}, and Galkowski--Zelditch studied lower bounds for Cauchy data on negatively curved surfaces \cite{GZ21}. Here the joint equations provide invariance under all commuting Hamiltonian flows, and isolation of a tangential branch prevents cancellation.

The construction in Theorem~\ref{thm:counterexample} concerns a single fixed curve at infinitely many eigenvalues of a single metric. Prescribed-nodal-set results for individual eigenfunctions, such as \cite{EPS}, address a interesting yet different question. Common curved hypersurfaces in higher-dimensional flat tori already arise from products \cite[\S1.2]{BR11}; our result concerns a closed surface. The underlying Liouville framework is standard \cite{GT}, and the exact construction uses elementary Floquet theory. The geometric conclusion is the failure of geodesic rigidity even for analytic metrics without Killing fields arbitrarily close to a flat metric.

Sections~\ref{sec:cauchy} and \ref{sec:connected-proof} prove the local Cauchy estimate and Theorem~\ref{thm:connected}. Sections~\ref{sec:nodal} and \ref{sec:sphere} retain the nodal-intersection application and spherical examples. Sections~\ref{sec:finite} and \ref{sec:counterexample} contain the finite-component extension and rigidity consequences, and the the exact nongeodesic counterexample and its interpretation.

\subsection*{Disclaimer on AI usage}
This paper uses generative AI models. Precisely, we used ChatGPT 6 Astra, Deepseek V4.1-Flash for the mathematical part and the writing, checking and graphing. Some other AI models are used in a nonessential way during the writing. Of course, the authors are responsible for the correctness of the whole content of this paper.

\subsection*{Acknowledgements}
The first author is supported by the research project funded by Wuhan University of Technology (Grant No.~202601hx0314).

Both authors would like to thank Prof. John Toth and Prof. Long Jin for helpful discussions relating to this work.

%% file: sections/core.tex
\section{A microlocal Cauchy estimate}\label{sec:cauchy}

We first prove a local estimate which does not require the joint fiber to be
connected. Its constants are uniform for the joint energy in a sufficiently
small fixed neighborhood of $E$. In particular, convergence $E(h)\to E$ is
sufficient; no estimate on its rate is needed.

We use two consequences of joint ellipticity. For joint eigenfunctions with
energy $e$ in a fixed bounded neighborhood of $E$, elliptic regularity gives
uniform bounds in every semiclassical Sobolev space and a fixed compact
microlocalization:
\begin{equation}\label{eq:core-compact}
 u=\Pi_hu+O(h^\infty)\norm{u}_{L^2(M)}
 \quad\text{in every }H_h^m(M),
 \qquad \Pi_h\in\Psi_h^{\mathrm{comp}}(M).
\end{equation}
Moreover, a compact microlocal cutoff supported where $|p-e|\ge c>0$
annihilates $u$ to $O(h^\infty)$, uniformly when the separation constant is
uniform. These assertions follow by applying an elliptic parametrix to
$\sum_k(P_k-e_k)^2$; see, for example, \cite{Zwo12}. Here and below an
$O(h^\infty)$ remainder has all the Sobolev bounds required for restriction.
For a compactly microlocalized operator $G_h$ we will use
\begin{equation}\label{eq:core-trace}
 \norm{R_HG_hw}_{L^2(H)}\le Ch^{-1/2}\norm{w}_{L^2(M)},
 \qquad R_Hw=w|_H.
\end{equation}
In local coordinates this follows from Cauchy--Schwarz in the bounded
semiclassical normal-frequency variable.

Choose an isolated transverse point $z_0\in\Lambda_E\cap T_H^*M$ as in
the hypothesis of Theorem~\ref{thm:connected}. Write $x=(t,x')$, with
$H=\{t=0\}$, and choose $j$ so that $\partial_{\xi_t}p_j\ne0$ near
$z_0$. Shrinking the tangential window $W$ and the selected neighborhood
if necessary, we may arrange that their relevant closures are compact.
The isolation hypothesis then persists for all $e$ sufficiently close to
$E$: every point of
\begin{equation}\label{eq:core-isolation}
 p^{-1}(e)\cap T_H^*M\cap r_H^{-1}(\overline{W_0}),
 \qquad W_0\Subset W,
\end{equation}
lies in the selected noncharacteristic neighborhood. Indeed, a contrary
sequence has a convergent subsequence by joint ellipticity and would
produce a competing point of the full fiber at $E$. This argument uses
the isolation of the \emph{entire} joint fiber, not just one of its
connected components.

All coordinate calculations below are made in Lebesgue measure. If the
ambient density is $m(t,x')\,dt\,dx'$, conjugate locally by $m^{1/2}$.
The conjugated operators have the same principal symbols and uniformly
bounded subprincipal symbols. The restriction of $m^{1/2}u$ differs from
the natural density realization of $u|_H$ by a smooth positive
multiplication operator. Thus the coordinate estimates give the stated
estimates in the original norms. Spatial cutoffs are taken equal to one
on a larger neighborhood of the chosen flow box; their errors are
microlocally negligible in the smaller neighborhoods used below.

\subsection{Factorization at the actual energy}

The following noncharacteristic reduction is the usual first-order
evolution factorization; compare \cite{BGT04,KTZ07,Zwo12}. We include its
subprincipal construction because it is important for the size of the
evolution error.

\begin{lemma}[Uniform parameter-dependent factorization]\label{lem:core-factor}
For $e$ in a sufficiently small neighborhood of $E$, there exist
$Q_{e,h}\in\Psi_h^0$, elliptic in the selected neighborhood, and a
tangential family $A_{e,h}(t)\in\Psi_h^1$ such that
\begin{equation}\label{eq:core-factor}
 Q_{e,h}(P_j-e_j)=hD_t-A_{e,h}(t)+h^2S_{e,h}
\end{equation}
microlocally there. Here $D_t=-i\partial_t$ and, after the compact
microlocal cutoffs, $S_{e,h}$ is uniformly bounded on $L^2$; it may be
taken in $\Psi_h^{-1}$. The tangential symbol is
\[
 a_{e,h}=a_e+h a_{1,e,h},\qquad a_e\text{ real},
 \qquad a_{1,e,h}\in S^0\text{ uniformly},
\]
and
\begin{equation}\label{eq:core-skew}
 A_{e,h}-A_{e,h}^*\in h\Psi_h^0
\end{equation}
uniformly in $e$ and $h$.
\end{lemma}

\begin{proof}
The implicit function theorem, applied with the parameter $e_j$, gives
\[
 r_e:=p_j-e_j=e_{0,e}(t,x',\xi_t,\xi')
          \bigl(\xi_t-a_e(t,x',\xi')\bigr),
 \qquad |e_{0,e}|\ge c>0.
\]
All symbols and their derivatives are uniform on fixed smaller
neighborhoods. Write the full Weyl symbol of $P_j-e_j$ as
$r_e+h p_{j,1,h}$ and set $q_{0,e}=e_{0,e}^{-1}$. With the convention
$H_a b=\{a,b\}$, the first Weyl correction is
\[
 d_{e,h}=q_{0,e}p_{j,1,h}+\frac{1}{2i}\{q_{0,e},r_e\}.
\]
Define
\[
 a_{1,e,h}(t,x',\xi')=-d_{e,h}(t,x',a_e(t,x',\xi'),\xi').
\]
The function $d_{e,h}+a_{1,e,h}$ vanishes on $\xi_t=a_e$, so smooth
division gives
\[
 d_{e,h}+a_{1,e,h}=(\xi_t-a_e)c_{e,h}.
\]
Taking $q_{1,e,h}=-e_{0,e}^{-1}c_{e,h}$, the Weyl product satisfies
\[
 (q_{0,e}+hq_{1,e,h})\#(r_e+hp_{j,1,h})
   =\xi_t-a_e-ha_{1,e,h}+O(h^2).
\]
This proves \eqref{eq:core-factor}. Extend $a_e$ as a real tangential
symbol, using fixed cutoffs outside the working neighborhood. The real
principal symbol and the bounded subprincipal term give
\eqref{eq:core-skew}. In particular, the construction factors
$p_j-e_j$ itself: the quantity $(e_j-E_j)/h$ never occurs.
\end{proof}

\subsection{Propagation and comparison with Dirichlet data}

\begin{proposition}[Local Cauchy estimate]\label{prop:local-cauchy}
At an isolated transverse point $z_0$ there are a fixed
$\chi\in C_c^\infty(T^*M)$, a uniformly bounded tangential operator
$B_H\in\Psi_h^0(H)$, a neighborhood $\mathcal V$ of $E$, and constants
$C,h_0>0$ such that
\begin{equation}\label{eq:core-local}
 \norm{\Op(\chi)u}_{L^2(M)}
 \le C\norm{B_H(u|_H)}_{L^2(H)}
       +Ch^{1/2}\norm{u}_{L^2(M)}
 \le C'\norm{u|_H}_{L^2(H)}
       +Ch^{1/2}\norm{u}_{L^2(M)}
\end{equation}
for every exact joint eigenfunction $P_ku=e_ku$, $e\in\mathcal V$, and
$0<h<h_0$. The principal symbol $\beta$ of $B_H$ is supported in $W_0$
and equals one near $r_H(z_0)$. The symbol $\chi$ equals one on a
nonempty relatively open flow box in $\Lambda_E$ through $z_0$.
\end{proposition}

\begin{proof}
Write $P=P_j$ and $p=p_j$ in this proof, and initially take
$\norm{u}_{L^2(M)}=1$. Choose $\beta\in C_c^\infty(W_0)$, equal to
one near $r_H(z_0)$. On the transversal $\{t=0\}$ choose a compact
cutoff
\[
 \alpha(x',\xi_t,\xi')=\beta(x',\xi')\vartheta(x',\xi_t,\xi')
\]
in the selected neighborhood, where $\vartheta=1$ near the graph
$\xi_t=a_e(0,x',\xi')$ over $\supp\beta$ for all $e\in\mathcal V$.
Extend $\alpha$ along the $H_p$ flow. Multiplying by a longitudinal
cutoff outside a larger fixed slab produces
$b\in C_c^\infty(T^*M)$ such that
\begin{equation}\label{eq:core-invariant}
 H_pb=0\quad\text{for }|t|<3\delta,
 \qquad b=1\quad\text{near }z_0.
\end{equation}
All supports lie inside the factorization neighborhood in this slab.
The number $\delta>0$ is fixed independently of $h$ and $e$.

Set $B_h=\Op(b)$ and $v=B_hu$. The commutator formula and
\eqref{eq:core-invariant} give
\[
 (P-e_j)v=[P,B_h]u=O(h^2)\quad\text{in the working slab}.
\]
This is a localized statement: the principal commutator vanishes in
the slab, whereas the derivatives of the longitudinal cutoff are
outside it. Insert nested cutoffs in the normal variable before
applying Lemma~\ref{lem:core-factor}. Pseudolocality makes the separated
cutoff errors $O(h^\infty)$ in the inner slab, and hence
\begin{equation}\label{eq:core-evolution}
 (hD_t-A_{e,h}(t))v=f_h,
 \qquad \norm{f_h}_{L^2((-2\delta,2\delta)\times\R^{n-1})}
       \le Ch^2.
\end{equation}
The real tangential principal symbol and \eqref{eq:core-skew} give the
energy estimate
\[
 \norm{v(t)}_{L^2_{x'}}
 \le C\norm{v(0)}_{L^2_{x'}}
    +\frac{C}{h}\int_{\min(0,t)}^{\max(0,t)}
                 \norm{f_h(s)}_{L^2_{x'}}\,ds,
 \qquad |t|\le\delta.
\]
One may obtain this directly by differentiating $\norm{v(t)}^2$;
the skew-adjoint part contributes at most $C\norm{v(t)}^2$.
Consequently,
\begin{equation}\label{eq:core-slab}
 \norm{v}_{L^2(\{|t|<\delta\})}
       \le C\norm{R_HB_hu}_{L^2(H)}+Ch.
\end{equation}

We next compare the localized trace with the original Dirichlet data.
Define the tangential symbol
\[
 \beta_e(t,x',\xi')=b(t,x',a_e(t,x',\xi'),\xi'),
 \qquad \beta_e(0,x',\xi')=\beta(x',\xi'),
\]
and let $\mathcal B_{e,h}(t)=\operatorname{Op}_{h,x'}^w(\beta_e(t))$. Smooth division in
the selected neighborhood gives
\begin{equation}\label{eq:core-division}
 b-\beta_e=c_e(\xi_t-a_e).
\end{equation}
Choose a compact microlocal cutoff $\Theta_h$ supported in that
neighborhood and equal to the identity near the relevant support of
$b$ at $H$ and near every point of \eqref{eq:core-isolation} over
$\supp\beta$. Joint ellipticity and the absence of competing branches
give
\begin{equation}\label{eq:core-boundary-cutoff}
 R_HB_hu=R_HB_h\Theta_hu+O(h^\infty),\qquad
 \operatorname{Op}_{h,x'}^w(\beta)R_Hu
       =R_H\mathcal B_{e,h}\Theta_hu+O(h^\infty).
\end{equation}
For the second equality, only the tangentially localized restriction
is being discarded outside the selected neighborhood. To justify it,
first use \eqref{eq:core-compact}; the restriction relation then retains
precisely the ambient covectors over $\supp\beta$. Those on the joint
fiber lie where $\Theta_h$ is the identity, and the remaining compact
set is uniformly joint elliptic. This also proves the asserted
$O(h^\infty)$ after taking the trace.

Quantizing \eqref{eq:core-division}, with nested compact microlocal
cutoffs, gives
\begin{equation}\label{eq:core-boundary-division}
 (B_h-\mathcal B_{e,h})\Theta_h
   =C_{e,h}(hD_t-A_{e,h})\Theta_h+hF_{e,h}+O(h^\infty),
\end{equation}
where $C_{e,h}$ and $F_{e,h}$ are uniformly bounded and compactly
microlocalized. Importantly, the cutoff state need not solve the
homogeneous equation to order $h^2$. Instead,
\[
 (hD_t-A_{e,h})\Theta_hu
   =Q_{e,h}[P,\Theta_h]u+O(h^2)=O(h)
\]
on the support needed in \eqref{eq:core-boundary-division}.
The commutator is $O(h)$, which is sufficient. Applying
\eqref{eq:core-trace} to \eqref{eq:core-boundary-division} and then
\eqref{eq:core-boundary-cutoff} yields
\begin{equation}\label{eq:core-boundary}
 R_HB_hu=\operatorname{Op}_{h,x'}^w(\beta)R_Hu+O(h^{1/2})_{L^2(H)}.
\end{equation}
After undoing the density conjugation this is the same assertion with
a uniformly bounded tangential operator $B_H$ of principal symbol
$\beta$.

Finally choose $\eta\in C_c^\infty((-\delta,\delta))$ equal to one
near zero and set $\chi(t,x',\xi)=\eta(t)b(t,x',\xi)$. Then
\[
 \Op(\chi)=\eta B_h+O(h)_{L^2\to L^2}.
\]
Equations \eqref{eq:core-slab} and \eqref{eq:core-boundary} prove
\eqref{eq:core-local}. This last cutoff makes the estimate global on
$M$. The construction makes $\chi=1$ on a small flow box through
$z_0$. Homogeneity restores the factor $\norm{u}_{L^2(M)}$ in the
error term.
\end{proof}

\section{The connected-fiber restriction theorem}\label{sec:connected-proof}

We recall explicitly the invariant-measure argument, so that its use of
connectedness and its lack of a nonresonance assumption are clear.

\begin{lemma}[The joint semiclassical measure]\label{lem:core-measure}
Suppose $p^{-1}(E)=\Lambda_E$ is compact, connected, and regular. Every
semiclassical measure of normalized joint eigenfunctions with
$E(h)\to E$ is the normalized Liouville measure $\mu_E$ on
$\Lambda_E$. In particular, for each fixed compactly supported symbol
$a$,
\begin{equation}\label{eq:core-measure}
 \norm{\Op(a)u_h}_{L^2(M)}^2
       \longrightarrow\int_{\Lambda_E}|a|^2\,d\mu_E.
\end{equation}
\end{lemma}

\begin{proof}
Compact microlocalization \eqref{eq:core-compact} ensures that every
subsequence has a further subsequence with a semiclassical measure of
total mass one. The joint equations and $E(h)\to E$ put its support in
$\Lambda_E$. Self-adjointness and the commutator identity give
\[
 \int H_{p_j}a\,d\mu=0,\qquad j=1,\ldots,n.
\]
Thus $\mu$ is invariant under all commuting Hamiltonian flows. At a
regular point their vector fields are independent and span
$T\Lambda_E$. In Liouville--Arnold angle coordinates they generate all
translations of the connected torus. Their only invariant probability
measure is normalized Haar measure, namely $\mu_E$. This is the
connected-fiber form of the standard joint-measure description in
\cite{TZ03}. Every convergent subsequence has this same limit, proving
\eqref{eq:core-measure} for the original family.
\end{proof}

\begin{proof}[Proof of Theorem~\ref{thm:connected}]
Apply Proposition~\ref{prop:local-cauchy} at the isolated transverse
point. Transversality implies that a small relatively open patch
$\Gamma_0\subset\Lambda_E\cap T_H^*M$ through that point has a flow
tube
\[
 \mathcal T=\{\exp(tH_{p_j})(z):z\in\Gamma_0,\ |t|<t_0\}
\]
which is open in $\Lambda_E$. The cutoff in the proposition can be
chosen equal to one on a smaller such tube. Consequently,
\[
 c_0:=\int_{\Lambda_E}|\chi|^2\,d\mu_E>0.
\]
Lemma~\ref{lem:core-measure} gives
$\norm{\Op(\chi)u_h}_{L^2(M)}\ge\sqrt{c_0/2}$ for sufficiently
small $h$. On the other hand, \eqref{eq:core-local} bounds this norm
by $C\norm{u_h|_H}_{L^2(H)}+Ch^{1/2}$. Absorbing the last term proves
the asserted positive lower bound.
\end{proof}

\begin{corollary}[Conditional $L^p$ restriction bounds]\label{cor:core-lp}
In the setting of Theorem~\ref{thm:connected}, assume in addition that
$H$ is compact and an independent estimate gives
$\norm{u_h|_H}_{L^\infty(H)}\le C$ for the normalized family. Then,
for every $1\le p\le\infty$,
\[
 c_p\le\norm{u_h|_H}_{L^p(H)}\le C_p.
\]
\end{corollary}
\begin{proof}
The upper bounds follow from the assumed sup-norm estimate. For
$1\le p\le2$ use
\[
 \norm{u_h|_H}_{L^2(H)}^2
 \le\norm{u_h|_H}_{L^\infty(H)}^{2-p}
       \norm{u_h|_H}_{L^p(H)}^p.
\]
For $2\le p\le\infty$, H\"older's inequality and the $L^2$ lower bound
give the result.
\end{proof}

\begin{remark}[Connectedness and nonresonance]\label{rem:core-nonresonance}
No arithmetic nonresonance assumption is used: invariance under all
$n$ flows already forces Haar measure on a connected regular fiber.
For a disconnected regular fiber, invariant probability measures can
distribute their mass among its torus components. The local estimate
above remains valid on an isolated sheet; the additional issue is
whether every possible component carrying mass is detected by the
restriction.
\end{remark}

\begin{remark}[Turning boundaries]\label{rem:core-turning}
The transverse hypothesis cannot hold at a smooth one-sided boundary
of the projected torus. Indeed, suppose locally
$\pi(\Lambda_E)\subset\{\rho\ge0\}$, with $H=\{\rho=0\}$ and
$d\rho\ne0$. At $z\in\Lambda_E\cap T_H^*M$, every $H_{p_j}$ flow
stays on $\Lambda_E$. The function
$t\mapsto\rho(\pi(\exp(tH_{p_j})z))$ has a local minimum at zero,
and hence $H_{p_j}\rho(z)=0$ for every $j$. Thus a lower-bound problem
at such a turning hypersurface requires an argument other than the
noncharacteristic Cauchy estimate used here.
\end{remark}

%% file: sections/applications.tex
\section{Sign changes from restriction lower bounds}\label{sec:nodal}

The restriction estimate has a nodal consequence when combined with decay of
period integrals.  For a complex-valued eigenfunction, it is preferable to count
sign changes of its real and imaginary parts separately: the union of their
zero sets is the entire curve if either part vanishes identically.

For a continuous real-valued function $f$ on a parametrized compact arc $I$,
let $N_{\mathrm{sc}}(f;I)$ be the supremum of the integers $N\geq 0$ for which
there are parameters $s_0<\cdots<s_N$ in $I$ satisfying
$f(s_i)f(s_{i+1})<0$ for $0\leq i<N$.  We set
$N_{\mathrm{sc}}(0;I)=0$ and allow the value $+\infty$.
For functions with isolated zeros this is the usual number of sign changes.

Whenever a square root of the Laplacian is used below, its symbol is understood
near the positive joint energy. A smooth spectral modification at low energy
puts it in the global semiclassical calculus without changing the eigenfunctions
under consideration.

\begin{proposition}[A sign-change criterion]\label{prop:sign-change-criterion}
Let $\gamma$ be a fixed compact smooth arc, and let $v_h$ be continuous
complex-valued functions on $\gamma$.  Suppose that
\begin{equation}\label{eq:nodal-sup-assumption}
 \sup_{0<h<h_0}\norm{v_h}_{L^\infty(\gamma)}<\infty,
\end{equation}
and that every nontrivial compact subarc $I$ in the interior of $\gamma$
satisfies
\begin{equation}\label{eq:nodal-local-assumptions}
 \liminf_{h\to0}\int_I |v_h|^2\,ds>0,
 \qquad
 \int_I v_h\,ds\longrightarrow0.
\end{equation}
Then
\begin{equation}\label{eq:nodal-sign-growth}
 N_{\mathrm{sc}}(\Repart v_h;\gamma)
 +N_{\mathrm{sc}}(\Impart v_h;\gamma)\longrightarrow\infty.
\end{equation}
In particular, if the functions $v_h$ are real-valued, their number of sign
changes tends to infinity.
\end{proposition}

\begin{proof}
Fix an integer $N\geq1$ and choose $N$ pairwise disjoint nontrivial compact
subarcs $I_1,\ldots,I_N$ in the interior of $\gamma$.
The first condition in \eqref{eq:nodal-local-assumptions} and
\eqref{eq:nodal-sup-assumption} imply, for each $j$, that
\begin{equation}\label{eq:nodal-lone}
 \int_{I_j} \bigl(|\Repart v_h|+|\Impart v_h|\bigr)\,ds
 \geq \int_{I_j}|v_h|\,ds
 \geq \frac{\int_{I_j}|v_h|^2\,ds}
 {\norm{v_h}_{L^\infty(\gamma)}}\geq c_j>0
\end{equation}
for all sufficiently small $h$.  If neither real nor imaginary part takes
both positive and negative values on $I_j$, then each has a constant weak
sign there.  Consequently,
\[
 \int_{I_j}\bigl(|\Repart v_h|+|\Impart v_h|\bigr)\,ds
 =\left|\Repart\int_{I_j}v_h\,ds\right|
  +\left|\Impart\int_{I_j}v_h\,ds\right|
 \leq\sqrt2\left|\int_{I_j}v_h\,ds\right|=o(1),
\]
contradicting \eqref{eq:nodal-lone}.  Thus each $I_j$ contributes a sign
transition of at least one of the two parts.  Since the arcs are disjoint,
the sum in \eqref{eq:nodal-sign-growth} is at least $N$ for all sufficiently
small $h$.  The integer $N$ was arbitrary.
\end{proof}

We next give hypotheses under which the three estimates in the proposition
follow from results on joint eigenfunctions.  In this paragraph $M$ is a
surface, $P_1(h)=\sqrt{-h^2\Delta_g}$, and the limiting first energy is $E_1=1$.
For a geodesic segment $\gamma$ parametrized by arclength $\tau$, set
\[
 \mathcal C_\gamma=\{(x,\xi):x\in\gamma,\ |\xi|_g^2=1\}.
\]
The admissibility condition of Wang--Wu \cite[Definition~1]{WW}, applied to
the moment map $(p_1^2,p_2)$, is
\begin{equation}\label{eq:ww-admissibility}
 |\partial_\tau p_2(x,\xi)|\geq c_\epsilon>0
 \quad\text{on}\quad
 \mathcal C_\gamma\cap\{|p_2-E_2|<\epsilon\}
\end{equation}
for some $\epsilon>0$.  Here $\partial_\tau$ differentiates at fixed cotangent
components in Fermi coordinates along $\gamma$, with $\tau$ the arclength
coordinate.  This condition is a separate assumption; it does not follow from
the noncharacteristic crossing condition of Theorem~\ref{thm:connected}.

\begin{corollary}[Nodal intersections on admissible geodesic segments]
\label{cor:nodal-geodesic}
Assume the connected regular joint-fiber and joint-eigenvalue hypotheses of
Theorem~\ref{thm:connected}, with $M$ a surface,
$P_1(h)=\sqrt{-h^2\Delta_g}$, $E_1=1$, and $E(h)=E+O(h)$.  Let $\gamma$ be a compact embedded
geodesic segment.  Suppose that
\begin{enumerate}
 \item every point of $\gamma$ belongs to $\pi(\Lambda_E)$;
 \item $d_\xi p_1\wedge d_\xi p_2\ne0$ on
 $\Lambda_E\cap T^*_\gamma M$;
 \item $r_\gamma:\Lambda_E\cap T^*_\gamma M\longrightarrow T^*\gamma$
 is injective;
 \item $\gamma$ satisfies \eqref{eq:ww-admissibility}.
\end{enumerate}
Then every normalized joint eigenfunction family in that theorem satisfies
\[
 N_{\mathrm{sc}}(\Repart u_h|_\gamma;\gamma)
 +N_{\mathrm{sc}}(\Impart u_h|_\gamma;\gamma)\longrightarrow\infty.
\]
For real-valued joint eigenfunctions this is a conclusion about the sign
changes of $u_h|_\gamma$ itself.
\end{corollary}

\begin{proof}
At any point $(x,\xi)\in\Lambda_E\cap T^*_\gamma M$, the two vectors
$\partial_\xi p_1$ and $\partial_\xi p_2$ span $T_xM$ by the fiber-rank
assumption.  Therefore at least one of the two Hamiltonian vector fields
crosses $\gamma$ at this point.  The first and third hypotheses allow us to
apply Theorem~\ref{thm:connected} to every nontrivial compact subarc
$I\subset\gamma^\circ$.  We obtain
\[
 \liminf_{h\to0}\norm{u_h}_{L^2(I)}>0.
\]
The fiber-rank hypothesis also gives
\begin{equation}\label{eq:nodal-local-sup}
 \norm{u_h}_{L^\infty(\gamma)}\leq C.
\end{equation}
Indeed, after joint elliptic localization to a neighborhood of the joint
fiber, the rank-one pointwise estimate on the energy circle gives this bound
on a neighborhood of each point of $\gamma$; see
\cite[Theorem~1 and the second assertion of Theorem~2]{WXb}.
Compactness of $\gamma$ makes the constant uniform.  The rank condition here
is imposed on all characteristic lifts over $\gamma$, rather than merely on
one of its branches.

Every fixed subsegment $I$ inherits \eqref{eq:ww-admissibility}.  Since
$P_1(h)^2=-h^2\Delta_g$ and $E_1(h)^2=1+O(h)$,
\cite[Theorem~2]{WW} applies and yields
\begin{equation}\label{eq:ww-fixed-subsegment}
 \left|\int_I u_h\,ds\right|\leq C_I h^{1/2}.
\end{equation}
Proposition~\ref{prop:sign-change-criterion} now applies to $u_h|_\gamma$.
\end{proof}

\begin{remark}
This proof uses unweighted periods over finitely many fixed subsegments, as
stated in \cite{WW}; it does not require a weighted extension of that theorem.
The constants $C_I$ in \eqref{eq:ww-fixed-subsegment} are allowed to depend on
$I$.  Consequently the argument gives qualitative divergence and does not,
by itself, give a uniform power-law lower bound for the number of sign
changes.  This is the same comparison between absolute mass and signed
periods that underlies nodal arguments such as \cite{JZ16}.
\end{remark}

\section{Spherical harmonics and the necessity of branch isolation}
\label{sec:sphere}

The sphere illustrates the connected-fiber restriction theorem and shows
why crossing alone is insufficient, even when the full joint fiber is
connected.  Equip $S^2$ with the round metric
\[
 g=d\theta^2+\sin^2\theta\,d\varphi^2,
 \qquad 0<\theta<\pi,
\]
and consider
\[
 P_1(h)=\sqrt{-h^2\Delta_{S^2}},\qquad
 P_2(h)=-ih\partial_\varphi.
\]
These operators commute, and their principal symbols are
\[
 p_1=\left(\xi_\theta^2+\frac{\xi_\varphi^2}{\sin^2\theta}\right)^{1/2},
 \qquad p_2=\xi_\varphi.
\]
The normalized complex spherical harmonics
$Y_\ell^m(\theta,\varphi)=F_{\ell,m}(\theta)e^{im\varphi}$ are joint
eigenfunctions, where $F_{\ell,m}$ is a real multiple of
$P_\ell^m(\cos\theta)$.

Fix $0<\mu<1$, put $h_\ell=(\ell+\tfrac12)^{-1}$, and choose an integer
$m_\ell$ nearest to $\mu/h_\ell$.  Then, for all sufficiently large $\ell$,
$0\leq m_\ell\leq\ell$ and
\begin{equation}\label{eq:sphere-energies}
 E_{1,\ell}=h_\ell\sqrt{\ell(\ell+1)}
 =\sqrt{1-h_\ell^2/4}=1+O(h_\ell^2),
 \qquad
 E_{2,\ell}=h_\ell m_\ell=\mu+O(h_\ell).
\end{equation}
Thus the joint eigenvalues obey the $O(h)$ energy-window condition.
The limiting joint fiber is
\begin{equation}\label{eq:sphere-joint-torus}
 \Lambda_\mu=\left\{|\xi|_g=1,\quad\xi_\varphi=\mu\right\},
 \qquad
 \xi_\theta^2+\frac{\mu^2}{\sin^2\theta}=1.
\end{equation}
It is a compact connected regular Lagrangian torus: the two signs of
$\xi_\theta$ join at the two turning latitudes.  Its projection is
\[
 \pi(\Lambda_\mu)=
 \{\arcsin\mu\leq\theta\leq\pi-\arcsin\mu\}.
\]

\subsection{Meridian restriction}
Let $I$ be a nontrivial compact interval strictly inside
$(\arcsin\mu,\pi-\arcsin\mu)$, and let
$H_I=\{(\theta,0):\theta\in I\}$.  For the local defining function
$\rho=\varphi$, one has
\[
 H_{p_2}\rho=1,
 \qquad
 r_{H_I}(\theta,0,\xi_\theta,\mu)=(\theta,\xi_\theta).
\]
The latter map is injective on the joint fiber.  Therefore
Theorem~\ref{thm:connected} gives
\begin{equation}\label{eq:sphere-meridian-lower}
 \norm{Y_\ell^{m_\ell}|_{H_I}}_{L^2(H_I)}\geq c_I>0
\end{equation}
for all sufficiently large $\ell$.  The fiber-rank condition also holds
over $H_I$, so the local pointwise bound is uniform there.  These are
nonturning-region estimates; the different behavior at the turning
latitudes is described, for example, in \cite{BDWZ12}.

For completeness, the real meridian restriction has order $h_\ell^{-1}$
sign changes on $I$.  Set
$w_\ell(\theta)=\sqrt{\sin\theta}\,F_{\ell,m_\ell}(\theta)$.  The associated
Legendre equation gives the exact equation
\begin{equation}\label{eq:sphere-sturm}
 w_\ell''+\frac1{h_\ell^2}
 \left(1-\frac{(h_\ell m_\ell)^2-h_\ell^2/4}{\sin^2\theta}\right)w_\ell=0.
\end{equation}
The coefficient in parentheses is bounded above and below by positive
constants on $I$ for all large $\ell$.  Sturm comparison with constant
positive coefficients therefore bounds the number of zeros in $I$ above
and below by constant multiples of $h_\ell^{-1}$, up to additive constants.
Every such zero is simple, by uniqueness for the second-order equation, and
is consequently a sign change.  This conclusion follows directly from the
separated equation and does not require admissibility in the sense of
\eqref{eq:ww-admissibility}.

\subsection{Cancellation on the equator}
Consider instead $H_{\mathrm{eq}}=\{\theta=\pi/2\}$ and
$\rho=\theta-\pi/2$.  On the same connected regular joint torus,
\[
 H_{p_1}\rho=\frac{\xi_\theta}{|\xi|_g}
 =\pm\sqrt{1-\mu^2}\ne0.
\]
Thus a commuting Hamiltonian already crosses the equator at every
characteristic lift.  If crossing by the second operator is desired, one
can equally use
\[
 \widetilde P_2(h)=P_2(h)+2P_1(h),\qquad
 H_{\widetilde p_2}\rho=\pm2\sqrt{1-\mu^2}\ne0.
\]
This replacement preserves the joint fiber, with limiting energy
$(1,\mu+2)$.

Now choose $m_\ell$ to be an integer nearest to $\mu/h_\ell$ among the
integers satisfying $\ell+m_\ell$ odd.  The distance to that parity class
is at most one, so \eqref{eq:sphere-energies} remains valid.  The identity
\[
 P_\ell^m(-t)=(-1)^{\ell+m}P_\ell^m(t)
\]
implies
\begin{equation}\label{eq:sphere-equator-zero}
 Y_\ell^{m_\ell}(\pi/2,\varphi)=0,
 \qquad
 \norm{Y_\ell^{m_\ell}|_{H_{\mathrm{eq}}}}_{L^2(H_{\mathrm{eq}})}=0.
\end{equation}
The failure of branch isolation is explicit:
\[
 r_{H_{\mathrm{eq}}}
 \left(\frac\pi2,\varphi,\pm\sqrt{1-\mu^2},\mu\right)
 =(\varphi,\mu).
\]
Both normal branches have the same tangential covector, and their
restrictions cancel exactly.  Hence a noncharacteristic crossing cannot
replace branch isolation in Theorem~\ref{thm:connected}.  This obstruction
occurs on a single connected regular torus, independently of the
additional issue of distributing mass among several components of a joint
fiber.

%% file: sections/finite_rigidity.tex
\section{Finite joint fibers and rigidity of persistent nodal curves}\label{sec:finite}

The local Cauchy estimate also applies when the full regular joint fiber has several connected components. This extension is useful for Liouville metrics, whose rotating joint fibers consist of four tori. We first prove the extension stated in Theorem~\ref{thm:finite}, and then use it to classify the curves on which exact joint eigenfunctions can vanish persistently. The lower bounds in this section follow from the microlocal argument of the paper; separation of variables will be used only for a uniform upper bound.

\subsection{Detecting every component of a regular fiber}

\begin{proof}[Proof of Theorem~\ref{thm:finite}]
Write the compact regular fiber as
\[
 \Lambda_E=\bigsqcup_{\nu=1}^{N}\Lambda_\nu.
\]
Each component is a compact Liouville torus. Let $\mu_\nu$ be its normalized Liouville measure. Joint ellipticity and the normalization of $u_h$ imply that every semiclassical defect measure is a probability measure supported on $\Lambda_E$. Invariance under all the commuting Hamiltonian flows, which span the tangent space of each component, gives
\begin{equation}\label{eq:finite-measure}
 \mu=\sum_{\nu=1}^{N}c_\nu\mu_\nu,
 \qquad c_\nu\ge0,\qquad \sum_{\nu=1}^{N}c_\nu=1.
\end{equation}
The weights need not be known and may depend on the subsequence.

For each $\nu$, choose the transverse point and tangential window supplied by the theorem. Shrink the window and the transverse patch so that their closures lie in the original neighborhoods. The isolation hypothesis concerns preimages in the \emph{entire} fiber: no other sheet, whether in $\Lambda_\nu$ or in another component, has a tangential covector in this smaller window. Applying Proposition~\ref{prop:local-cauchy}, with the transverse Hamiltonian chosen for this component, gives a cutoff $\chi_\nu$ in a small flow tube such that
\begin{equation}\label{eq:finite-cauchy}
 \norm{\Op(\chi_\nu)u_h}_{L^2(M)}
 \le C_\nu\norm{u_h|_H}_{L^2(H)}+o(1).
\end{equation}
We choose the support disjoint from the other components and the cutoff nonzero on an open subset of $\Lambda_\nu$. Consequently,
\[
 d_\nu:=\int_{\Lambda_\nu}|\chi_\nu|^2\,d\mu_\nu>0.
\]
All cutoffs and constants can be chosen for nearby joint energies, as in the local estimate. In particular, convergence $E(h)\to E$ suffices here.

If the asserted restriction lower bound failed, a sequence with $h\to0$ and $\norm{u_h|_H}_{L^2(H)}\to0$ would have a subsequence with defect measure \eqref{eq:finite-measure}. Equation \eqref{eq:finite-cauchy} would then imply, for every $\nu$,
\[
 0=\lim_{h\to0}\norm{\Op(\chi_\nu)u_h}_{L^2(M)}^2
   =c_\nu d_\nu.
\]
Since $d_\nu>0$, every $c_\nu$ would vanish, contradicting their sum being one. This proves the theorem without a condition selecting a component or prescribing its mass.
\end{proof}

\begin{corollary}[Uniformity in the joint energy]\label{cor:finite-uniform}
Suppose that the geometric hypotheses of Theorem~\ref{thm:finite} hold at $E$. There are a neighborhood $U$ of $E$ and constants $c,h_0>0$ such that every exact joint eigenfunction satisfying
\[
 P_j(h)u=E'_ju\quad(1\le j\le n),\qquad E'\in U,
 \qquad 0<h<h_0,
\]
obeys
\begin{equation}\label{eq:finite-uniform}
 \norm{u|_H}_{L^2(H)}\ge c\norm{u}_{L^2(M)}.
\end{equation}
The same conclusion holds with $H$ replaced by a fixed small patch containing the transverse points used in the hypothesis.
\end{corollary}

\begin{proof}
Otherwise, choosing successively smaller energy neighborhoods and smaller upper bounds on $h$, one could find normalized exact joint eigenfunctions with
\[
 h_k\to0,\qquad E'_k\to E,\qquad
 \norm{u_k|_H}_{L^2(H)}\to0.
\]
The proof just given applies to this sequence and yields a contradiction. For a patch, choose all the local Cauchy cutoffs over that patch. Uniformity of their estimates follows by factoring at the actual nearby energy: the transverse derivative remains bounded away from zero, and the tangential windows continue to exclude the other sheets after shrinking. No estimate of $E'_k-E$ by $O(h_k)$ is required.
\end{proof}

\begin{remark}
A transverse isolated patch on only one component does not suffice for a lower bound on arbitrary normalized joint eigenfunctions: that component can have weight zero in \eqref{eq:finite-measure}. Requiring such a patch on every component is what removes the need to determine the weights. The example in Theorem~\ref{thm:counterexample} will fail precisely at this point: the two components carrying the cancelling waves have the same tangential covectors along the common nodal curve.
\end{remark}

\subsection{The four branches of a Liouville metric}

Let $M=\T^2=(\R/2\pi\mathbb Z)^2$ and let
\begin{equation}\label{eq:liouville-metric}
 g=F(x,y)(dx^2+dy^2),\qquad F(x,y)=a(x)+b(y)>0,
\end{equation}
where $a,b$ are real smooth periodic functions. The commuting operators
\begin{equation}\label{eq:liouville-operators}
 L=-\Delta_g=-F^{-1}(\partial_x^2+\partial_y^2),
 \qquad Q=F^{-1}(-b(y)\partial_x^2+a(x)\partial_y^2)
\end{equation}
are self-adjoint on their natural domains in $L^2(M,F\,dx\,dy)$. Their semiclassical principal symbols are
\[
 \ell=\frac{\xi_x^2+\xi_y^2}{F},\qquad
 q=\frac{b\xi_x^2-a\xi_y^2}{F}.
\]
For $h>0$ and $\kappa\in\R$, define the exact joint eigenspace
\begin{equation}\label{eq:liouville-space}
 \E_{h,\kappa}
 =\{u\in C^\infty(M;\C):h^2Lu=u,\ h^2Qu=\kappa u\}.
\end{equation}
We work in the rotating interval
\begin{equation}\label{eq:liouville-rotating}
 \mathcal R=(-\min a,\min b),
\end{equation}
which is nonempty because $\min a+\min b>0$. For $\kappa\in\mathcal R$, put
\begin{equation}\label{eq:liouville-AB}
 A_\kappa(x)=\sqrt{a(x)+\kappa},\qquad
 B_\kappa(y)=\sqrt{b(y)-\kappa}.
\end{equation}
The full joint fiber $\{\ell=1,q=\kappa\}$ is the disjoint union of the four graph tori
\begin{equation}\label{eq:liouville-four-tori}
 \Lambda_{\sigma\tau}(\kappa)
 =\{(x,y;\sigma A_\kappa(x),\tau B_\kappa(y))\},
 \qquad \sigma,\tau\in\{+1,-1\}.
\end{equation}
They are regular, since
\[
 \det\frac{\partial(\ell,q)}{\partial(\xi_x,\xi_y)}
 =-\frac{4\xi_x\xi_y}{F}\ne0.
\]

For completeness, these second-order differential operators fit the first-order framework of the main theorem by setting
\begin{equation}\label{eq:liouville-first-order}
 P_1(h)=(1+h^2L)^{1/2},\qquad
 P_2(h)=h^2Q(1+h^2L)^{-1/2}.
\end{equation}
They commute and are self-adjoint, with real principal symbols
\[
 p_1=\sqrt{1+\ell},\qquad p_2=\frac{q}{\sqrt{1+\ell}}.
\]
The symbol map $(\ell,q)\mapsto(p_1,p_2)$ has invertible differential. The space \eqref{eq:liouville-space} is therefore exactly the joint eigenspace at energy $(\sqrt2,\kappa/\sqrt2)$ for \eqref{eq:liouville-first-order}, and its classical fiber is \eqref{eq:liouville-four-tori}.

Let $\gamma$ be a smooth regular curve. A point of $\gamma$ with nonzero tangent $v=(v_x,v_y)$ is called \emph{separating at energy $\kappa$} if
\begin{equation}\label{eq:separation}
 v_xv_y\ne0,\qquad
 A_\kappa(x)^2v_x^2\ne B_\kappa(y)^2v_y^2.
\end{equation}
This is independent of the regular parametrization. Indeed, the four tangential covectors in \eqref{eq:liouville-four-tori} are represented by
\begin{equation}\label{eq:liouville-tangential}
 \sigma A_\kappa v_x+\tau B_\kappa v_y,
 \qquad \sigma,\tau\in\{+1,-1\},
\end{equation}
and they are pairwise distinct exactly when \eqref{eq:separation} holds.

\begin{lemma}[Restriction near a separating point]\label{lem:liouville-restriction}
Suppose that \eqref{eq:separation} holds at a point of $\gamma$ at energy $\kappa_0\in\mathcal R$. There are a fixed compact subarc $\gamma_0\subset\gamma$, an open interval $J\Subset\mathcal R$ containing $\kappa_0$, and constants $c,h_0>0$ such that
\begin{equation}\label{eq:liouville-local}
 \norm{u|_{\gamma_0}}_{L^2(\gamma_0)}
 \ge c\norm{u}_{L^2(M)},
 \qquad u\in\E_{h,\kappa},\quad \kappa\in J,\quad 0<h<h_0.
\end{equation}
Curve norms are taken with respect to metric arclength.
\end{lemma}

\begin{proof}
At the specified point, choose four disjoint neighborhoods of the tangential covectors \eqref{eq:liouville-tangential}. Shrinking the base arc makes each neighborhood isolate its sheet from the entire four-component fiber. On each branch, the vectors $\partial_\xi p_1,\partial_\xi p_2$ span the base tangent space. If $\rho$ is a local defining function for the curve, at least one of $H_{p_1}\rho,H_{p_2}\rho$ is therefore nonzero. The transverse Hamiltonian may depend on the branch, which is permitted in Theorem~\ref{thm:finite}. Thus every component has a transverse isolated patch over the same short arc. Apply Corollary~\ref{cor:finite-uniform} to the system \eqref{eq:liouville-first-order}. This proves the estimate for all nearby $\kappa$.
\end{proof}

We will also use the following elementary upper bound. Its proof is included to keep the subsequent stability assertion independent of additional pointwise estimates.

\begin{lemma}[Uniform upper bound in a rotating window]\label{lem:liouville-sup}
For every compact $K\Subset\mathcal R$, there are $C_K,h_0>0$ such that
\begin{equation}\label{eq:liouville-sup}
 \norm{u}_{L^\infty(M)}\le C_K\norm{u}_{L^2(M)},
 \qquad u\in\E_{h,\kappa},\quad \kappa\in K,\quad 0<h<h_0.
\end{equation}
\end{lemma}

\begin{proof}
The identities $-\partial_x^2=aL+Q$ and $-\partial_y^2=bL-Q$ give
\begin{equation}\label{eq:liouville-separated}
 -h^2u_{xx}=A_\kappa^2u,\qquad
 -h^2u_{yy}=B_\kappa^2u.
\end{equation}
On a fixed coordinate period, each scalar equation has exact independent solutions with uniform expansions
\[
 f_\sigma(x)=A_\kappa(x)^{-1/2}e^{i\sigma S_x(x)/h}+O(h),
 \qquad
 g_\tau(y)=B_\kappa(y)^{-1/2}e^{i\tau S_y(y)/h}+O(h),
\]
where $S_x'=A_\kappa$, $S_y'=B_\kappa$, and the errors are in $C^0$, uniformly for $\kappa\in K$. To verify this elementary WKB estimate, the substitution $t=S_x(x)$, $f=A_\kappa^{-1/2}w$ gives
\[
 w_{tt}+h^{-2}w=r_\kappa(t)w,\qquad
 r_\kappa=\frac{A_\kappa''}{2A_\kappa^3}
             -\frac{3(A_\kappa')^2}{4A_\kappa^4}.
\]
The solution with initial data $w(0)=1$, $w_t(0)=i\sigma/h$ satisfies
\[
 w(t)=e^{i\sigma t/h}
       +h\int_0^t\sin((t-s)/h)r_\kappa(s)w(s)\,ds.
\]
Boundedness of $r_\kappa$ and Gronwall's inequality give the asserted expansion; the $y$ equation is identical in form.

Every solution of both equations in \eqref{eq:liouville-separated} has an exact representation
\[
 u(x,y)=\sum_{\sigma,\tau=\pm1}c_{\sigma\tau}f_\sigma(x)g_\tau(y).
\]
Indeed, expand first in the $x$ basis; its two coefficient functions must solve the $y$ equation. The scalar basis functions themselves need not be periodic. Their Gram matrices on a full coordinate period have positive diagonal entries bounded away from zero and off-diagonal entries $O(h)$: integration by parts applies to the leading off-diagonal phases because $A_\kappa,B_\kappa$ are bounded away from zero. Hence the tensor-product Gram matrix, together with the positive bounds on $F$, gives
\[
 \sum_{\sigma,\tau}|c_{\sigma\tau}|^2
 \le C_K\norm{u}_{L^2(M)}^2.
\]
The basis functions are uniformly bounded, so this coefficient estimate proves \eqref{eq:liouville-sup}.
\end{proof}

\subsection{Persistent curves and stable restriction}

Define the closed nowhere-zero one-forms
\begin{equation}\label{eq:liouville-forms}
 \alpha_\kappa^\pm=A_\kappa(x)\,dx\pm B_\kappa(y)\,dy.
\end{equation}
A \emph{phase front} is a leaf of $\ker\alpha_\kappa^+$ or $\ker\alpha_\kappa^-$. Locally, it is a level curve of $S_x\pm S_y$; these primitives need not be single-valued on the torus. All curves below are fixed independently of $h$.

\begin{proposition}[Rigidity of persistent nodal curves]\label{prop:rigidity}
Let $\gamma$ be a connected smooth embedded regular curve. Suppose that normalized exact joint eigenfunctions $u_j\in\E_{h_j,\kappa_j}$ satisfy
\[
 h_j\to0,\qquad \kappa_j\to\kappa\in\mathcal R,
 \qquad u_j|_\gamma\equiv0.
\]
Then $\gamma$ is contained in a leaf of exactly one of the four foliations
\begin{equation}\label{eq:liouville-four-foliations}
 \ker dx,\qquad \ker dy,\qquad
 \ker\alpha_\kappa^+,\qquad\ker\alpha_\kappa^-.
\end{equation}
In particular, it is a coordinate curve or a phase front of the limiting joint energy. If $\gamma$ is closed, it is a closed leaf.
\end{proposition}

\begin{proof}
A separating point at the limiting energy would, by Lemma~\ref{lem:liouville-restriction}, give a fixed subarc on which the restriction norms of $u_j$ are bounded below for all large $j$. This contradicts their vanishing. Thus, at every point of a regular parametrization $\gamma(t)=(x(t),y(t))$,
\[
 x'y'=0\quad\text{or}\quad
 A_\kappa(x)^2x'^2=B_\kappa(y)^2y'^2.
\]
These are precisely the four tangent lines in \eqref{eq:liouville-four-foliations}. Since $A_\kappa,B_\kappa>0$, the four lines are distinct at every point. The continuous tangent line of $\gamma$ therefore selects a locally constant label, which is constant because $\gamma$ is connected. Uniqueness for integral curves of the resulting smooth line field places $\gamma$ in one leaf. A closed embedded curve tangent to that field is itself a closed leaf.
\end{proof}

\begin{corollary}[A contractible closed curve cannot persist]\label{cor:contractible}
Let $\gamma$ be a fixed smooth embedded contractible closed curve and let $K\Subset\mathcal R$ be compact. There are $c,h_0>0$ such that
\begin{equation}\label{eq:contractible}
 \norm{u|_\gamma}_{L^2(\gamma)}\ge c\norm{u}_{L^2(M)},
 \qquad u\in\E_{h,\kappa},\quad\kappa\in K,\quad0<h<h_0.
\end{equation}
In particular, $\gamma$ cannot be a persistent nodal component in this rotating window.
\end{corollary}

\begin{proof}
Fix $\kappa\in K$. If $\gamma$ had no separating point, the line-field argument in the preceding proof would make it a closed leaf of one of \eqref{eq:liouville-four-foliations}. Coordinate leaves are noncontractible. On a phase-front leaf,
\[
 A_\kappa(x)x'\pm B_\kappa(y)y'=0.
\]
Regularity implies that neither $x'$ nor $y'$ vanishes; both have constant sign. For a parametrization by $\R/\mathbb Z$ and a lift $(\widetilde x,\widetilde y)$, its winding numbers therefore satisfy
\[
 2\pi m=\int_0^1\widetilde x'(t)\,dt\ne0,
 \qquad
 2\pi n=\int_0^1\widetilde y'(t)\,dt\ne0.
\]
Such a leaf is also noncontractible, a contradiction. Hence a separating point exists at every $\kappa\in K$. Apply Lemma~\ref{lem:liouville-restriction} at these points and take a finite subcover of their energy neighborhoods. The minimum of the finitely many lower-bound constants and frequency thresholds gives \eqref{eq:contractible}.
\end{proof}

\begin{corollary}[Action resonance for closed noncoordinate curves]\label{cor:resonance}
Under the hypotheses of Proposition~\ref{prop:rigidity}, suppose that $\gamma$ is closed and noncoordinate. Define the positive actions
\begin{equation}\label{eq:liouville-actions}
 I_x(\kappa)=\int_0^{2\pi}A_\kappa(x)\,dx,
 \qquad
 I_y(\kappa)=\int_0^{2\pi}B_\kappa(y)\,dy.
\end{equation}
The winding numbers $(m,n)$ of $\gamma$ satisfy $mn\ne0$ and, for one fixed sign,
\begin{equation}\label{eq:liouville-resonance}
 mI_x(\kappa)\pm nI_y(\kappa)=0.
\end{equation}
Consequently $I_x(\kappa)/I_y(\kappa)=|n/m|\in\mathbb Q$. If this action ratio is irrational, every smooth embedded closed persistent nodal curve is a coordinate circle.
\end{corollary}

\begin{proof}
By Proposition~\ref{prop:rigidity}, the curve is a leaf of $\ker\alpha_\kappa^\pm$ for one sign. The preceding proof shows that $mn\ne0$. Integrating its annihilating one-form and using periodicity gives
\[
 0=\int_\gamma\alpha_\kappa^\pm
   =mI_x(\kappa)\pm nI_y(\kappa),
\]
which proves the assertions.
\end{proof}

\begin{remark}
The phase coordinates
\[
 X(x)=\int_0^x A_\kappa(t)\,dt,
 \qquad Y(y)=\int_0^y B_\kappa(t)\,dt
\]
identify the torus with a rectangular torus of side lengths $I_x(\kappa)$ and $I_y(\kappa)$. Its phase fronts become $X\pm Y=\mathrm{constant}$. Thus the resonance is exactly the condition for these leaves to close. It is necessary for persistence of a closed noncoordinate curve, but does not itself construct an exact eigenfunction vanishing on a chosen leaf.
\end{remark}

\begin{proposition}[Stable uniqueness from restriction]\label{prop:stable}
Let $\gamma$ be a fixed compact smooth embedded regular curve, possibly with boundary, and let $K\Subset\mathcal R$. Suppose that $\gamma$ has a separating point for every $\kappa\in K$. There are $c,C,h_0>0$ such that
\begin{equation}\label{eq:liouville-two-sided}
 c\norm{u}_{L^2(M)}\le\norm{u|_\gamma}_{L^2(\gamma)}
 \le C\norm{u}_{L^2(M)},
 \qquad u\in\E_{h,\kappa},\quad\kappa\in K,\quad0<h<h_0.
\end{equation}
Restriction is therefore injective on each such joint eigenspace, with a uniformly bounded inverse on its range. In particular, for $u,v$ in the same $\E_{h,\kappa}$,
\begin{equation}\label{eq:liouville-stability}
 \norm{u-v}_{L^2(M)}
 \le c^{-1}\norm{(u-v)|_\gamma}_{L^2(\gamma)}.
\end{equation}
\end{proposition}

\begin{proof}
The finite-cover argument using Lemma~\ref{lem:liouville-restriction} proves the lower estimate uniformly on $K$. Lemma~\ref{lem:liouville-sup} gives the upper estimate because
\[
 \norm{u|_\gamma}_{L^2(\gamma)}
 \le\operatorname{length}_g(\gamma)^{1/2}\norm{u}_{L^\infty(M)}.
\]
Finally, $u-v$ belongs to the same exact joint eigenspace, so the lower estimate applied to this difference gives \eqref{eq:liouville-stability}.
\end{proof}

These statements allow every vector in a joint eigenspace, including linear combinations when the joint eigenvalue has multiplicity. They concern exact joint eigenfunctions in a regular rotating window; no assertion is made for arbitrary Laplace eigenfunctions, or for differences of functions with different joint eigenvalues. The construction in the next section realizes one of the noncoordinate phase fronts permitted by Proposition~\ref{prop:rigidity}, and shows that this exceptional foliation need not consist of geodesics.

%% file: sections/counterexample.tex
\section{A nongeodesic common nodal circle}\label{sec:counterexample}

We now construct the exceptional noncoordinate phase front allowed by Section~\ref{sec:finite}. The existence argument is exact and independent of the restriction lower bound. The final calculation identifies the transverse branches whose traces cannot be separated.

\subsection{The Floquet construction}

Let $f>0$ be a nonconstant smooth function with period $\ell=2\pi/N$, where $N\ge3$ is an integer. Let $r,s$ be relatively prime positive integers with $r\ne s$. On $\T^2$, set
\begin{equation}\label{eq:general-metric}
 a(x)=r^2f(rx),\qquad b(y)=s^2f(sy),\qquad
 F=a+b,\qquad g=F(dx^2+dy^2),
\end{equation}
and define
\begin{equation}\label{eq:construction-operators}
 L=-\Delta_g=-F^{-1}(\partial_x^2+\partial_y^2),\qquad
 Q=F^{-1}(-b\partial_x^2+a\partial_y^2).
\end{equation}
The operators commute and are formally self-adjoint for $dV_g=F\,dx\,dy$.

\begin{proposition}[Matching separated equations]\label{prop:construction}
For the metric \eqref{eq:general-metric}, there exist real normalized joint eigenfunctions $u_j$ with
\[
 Lu_j=\lambda_j u_j,\qquad Qu_j=0,\qquad \lambda_j\to\infty,
\]
at distinct eigenvalues, for which
\begin{equation}\label{eq:general-curve}
 \gamma_{r,s}(t)=(st,rt)\pmod{2\pi}
\end{equation}
is a regular connected nodal component. This embedded circle has winding $(s,r)$ and signed geodesic curvature
\begin{equation}\label{eq:general-curvature}
 k_g(t)=\frac{(r^2-s^2)f'(rs t)}
 {2(r^2+s^2)f(rs t)^{3/2}}.
\end{equation}
It is not a geodesic, and $g$ has no nonzero Killing field. If $f$ is real analytic, so are the metric and the eigenfunctions.
\end{proposition}

\begin{proof}
\emph{The twisted one-dimensional problem.}
Put $\omega=e^{2\pi i/N}$ and consider the operator
\[
 T\phi=-f^{-1}\phi'',\qquad
 D(T)=\{\phi\in H^2([0,\ell]):
 \phi(\ell)=\omega\phi(0),\ \phi'(\ell)=\omega\phi'(0)\}
\]
on $L^2([0,\ell],f(t)\,dt)$. It is associated with the closed quadratic form $\int_0^\ell|\phi'|^2dt$ on the twisted $H^1$ domain. Since $|\omega|=1$, the boundary terms cancel; the operator is self-adjoint. Compactness of the form-domain embedding gives compact resolvent. The only constant function obeying the boundary condition is zero, so $T$ is strictly positive. Choose eigenpairs at distinct eigenvalues
\begin{equation}\label{eq:bloch}
 -\phi_j''=\lambda_j f\phi_j,\qquad \lambda_j\to\infty.
\end{equation}
Extend by $\phi_j(t+\ell)=\omega\phi_j(t)$. The matching values and derivatives, together with the differential equation, give a smooth global solution. It is $2\pi$-periodic because $\omega^N=1$.

The multiplier is nonreal, so $\phi_j$ and $\overline{\phi_j}$ are independent solutions: otherwise they would have the same multiplier, whereas $\omega\ne\overline\omega$. Their nonzero Wronskian gives a constant nonzero current
\begin{equation}\label{eq:current}
 J_j=\Impart\bigl(\phi_j'(t)\overline{\phi_j(t)}\bigr)\ne0.
\end{equation}
In particular, $\phi_j$ never vanishes.

\emph{Exact joint eigenfunctions and the common component.}
Define the real, doubly periodic function
\begin{equation}\label{eq:U}
 U_j(x,y)=\Impart\bigl(\phi_j(rx)\overline{\phi_j(sy)}\bigr).
\end{equation}
Equation \eqref{eq:bloch} yields
\[
 -\partial_x^2U_j=\lambda_j r^2f(rx)U_j=\lambda_j aU_j,
 \qquad
 -\partial_y^2U_j=\lambda_j s^2f(sy)U_j=\lambda_j bU_j.
\]
Consequently,
\begin{equation}\label{eq:exact-joint}
 LU_j=\lambda_jU_j,\qquad QU_j=0.
\end{equation}
Along \eqref{eq:general-curve}, both arguments in \eqref{eq:U} equal $rs t$, and hence
\begin{equation}\label{eq:regular-zero}
 U_j|_{\gamma_{r,s}}=0,\qquad
 dU_j|_{\gamma_{r,s}}=J_j(r\,dx-s\,dy)\ne0.
\end{equation}
Coprimality of $r,s$ makes \eqref{eq:general-curve} embedded. The implicit function theorem and compactness give a tubular neighborhood in which its nodal set is exactly this circle. Thus it is a regular connected component of the full nodal set. Normalize $U_j$ in $L^2(M,g)$ to obtain $u_j$. Distinct eigenvalues make the normalized functions orthogonal. Analyticity follows from the analytic differential equation when $f$ is analytic.

\emph{Geodesic curvature.}
Write $D=r^2+s^2$ and orient the curve by increasing $t$. Its Euclidean unit tangent and left normal are
\[
 T_0=(s,r)/\sqrt D,\qquad n_0=(-r,s)/\sqrt D.
\]
For a conformal metric $g=F(dx^2+dy^2)$, the signed geodesic curvature is
\[
 k_g=F^{-1/2}\bigl(k_0-\tfrac12\partial_{n_0}\log F\bigr).
\]
Here $k_0=0$, while on the curve
\[
 F=D f(rs t),\qquad
 \partial_{n_0}F=(s^2-r^2)\sqrt D\,f'(rs t).
\]
Substitution proves \eqref{eq:general-curvature}. Since $r\ne s$ and $f$ is nonconstant, this curvature is not identically zero.

\emph{Absence of Killing fields.}
If $X=P\partial_x+R\partial_y$ is a Killing field, the trace-free parts of $\mathcal L_Xg=0$ give
\[
 P_x=R_y,\qquad P_y=-R_x.
\]
Thus $P,R$ are harmonic periodic functions and hence constants. The remaining equation becomes
\[
 P r^3f'(rx)+R s^3f'(sy)=0\qquad\text{for all }(x,y).
\]
Integrating first in $y$ and then in $x$, and using that $f'$ is not identically zero, gives $P=R=0$.
\end{proof}

\begin{proof}[Proof of Theorem~\ref{thm:counterexample}]
Apply Proposition~\ref{prop:construction} with $N=3$, $r=1$, $s=2$, and $f(t)=2+\varepsilon\cos3t$. Positivity follows from $0<\varepsilon<2$. Formula \eqref{eq:general-curvature} becomes \eqref{eq:explicit-curvature}.
\end{proof}

\paragraph{Why the choices matter.}
The small period makes a nonreal Floquet multiplier into a periodic solution on the original circle. The factors $r^2,s^2$ exactly compensate for differentiating $\phi_j(rx)$ and $\phi_j(sy)$. Matching their arguments along $rx=sy$ produces cancellation for every $j$, while $r\ne s$ removes geodesicity. Finite translation symmetries remain, but there is no nonzero Killing field; the persistence mechanism is the matching of scalar equations. Since the common curve is not geodesic, it cannot be a curve fixed pointwise by a nontrivial isometric reflection.

\begin{remark}[A simpler example with a Killing field]\label{rem:warped}
The broad existence phenomenon also has a one-variable version. For $g_0=f(x)(dx^2+dy^2)$ and any $x_0$ with $f'(x_0)\ne0$, rephase each of the same Floquet solutions so that $\phi_j(x_0)>0$ is real. Then $v_j(x,y)=\Impart\phi_j(x)$ satisfies $-\Delta_{g_0}v_j=\lambda_jv_j$, and $x=x_0$ is a regular nodal circle because $v_j'(x_0)=J_j/\phi_j(x_0)\ne0$. Its geodesic curvature has magnitude $|f'(x_0)|/(2f(x_0)^{3/2})$. This direct comparison explains why the additional features of Proposition~\ref{prop:construction} matter: it removes the Killing field, makes the curve noncoordinate, and places the joint sequence in a rotating energy window.
\end{remark}

\input{sections/phase_interpretation}
\input{sections/quantum_limit}
\input{sections/counterexample_figures}

%% file: sections/phase_interpretation.tex
\subsection{Action resonance and transverse branch cancellation}

For \eqref{eq:general-metric}, the energy $\kappa=0$ lies in the rotating interval
$(-r^2\min f,s^2\min f)$. With $h_j=\lambda_j^{-1/2}$, the constructed
sequence lies in $\E_{h_j,0}$. Its full joint fiber has the four graph
components \eqref{eq:liouville-four-tori}, so the finite-component formulation
of Section~\ref{sec:finite} is the appropriate one.

Along $\gamma_{r,s}(t)=(st,rt)$, the one-form
\[
 \alpha_0^-=r\sqrt{f(rx)}\,dx-s\sqrt{f(sy)}\,dy
\]
annihilates the tangent. Thus the common nodal circle is a phase front of
the kind allowed by Proposition~\ref{prop:rigidity}. If
\[
 I=\int_0^{2\pi}\sqrt{f(t)}\,dt,
\]
periodicity and change of variables give
\[
 I_x(0)=rI,\qquad I_y(0)=sI,
 \qquad sI_x(0)-rI_y(0)=0.
\]
The winding $(s,r)$ therefore satisfies the resonance of
Corollary~\ref{cor:resonance}. Resonance closes the phase front; the
Floquet construction supplies exact eigenfunctions vanishing on it.

The failure of the restriction hypothesis is particularly transparent in
Theorem~\ref{thm:counterexample}. On $\gamma(t)=(2t,t)$, put
$A=\sqrt{2+\varepsilon\cos6t}$, so $B=2A$. The four covectors restricted
to the tangent $(2,1)$ give
\begin{equation}\label{eq:counterexample-momenta}
 (\sigma A,\tau B)(2,1)=2A(\sigma+\tau)
 \in\{4A,0,0,-4A\}.
\end{equation}
The mixed branches are distinct in $T^*M$, but both have zero tangential
covector. They cannot be separated by a pseudodifferential cutoff on the
curve. For the local defining function $\rho=x-2y$ and the Laplace symbol
$p_L=(\xi_x^2+\xi_y^2)/F_\varepsilon$, these same branches satisfy
\begin{equation}\label{eq:counterexample-crossing}
 H_{p_L}\rho\big|_{\Lambda_{+-}(0)}=\frac{2}{A},
 \qquad
 H_{p_L}\rho\big|_{\Lambda_{-+}(0)}=-\frac{2}{A}.
\end{equation}
Both flows cross the curve transversely. The obstruction is therefore
the coincidence of tangential covectors, despite nonzero geodesic curvature
and transverse Hamiltonian crossing. The two same-sign branches remain
individually observable; the next proposition uses the Cauchy estimate to
show that they carry no limiting mass in this sequence.

The distinction from Section~\ref{sec:sphere} should be kept explicit.
The spherical example violates branch isolation on one connected torus.
Here the cancellation occurs between two components of a four-component
regular fiber. It exhibits failure of the hypothesis that every component
be detectable in Theorem~\ref{thm:finite}.

%% file: sections/quantum_limit.tex
\subsection{The quantum limit selected by the common nodal curve}

The local Cauchy estimate determines the semiclassical distribution of the constructed eigenfunctions, although the full joint fiber is disconnected. This gives a direct application of the restriction method to the example, without requiring asymptotic information about its Floquet eigenfunctions.

\begin{proposition}[A unique limit on two components]\label{prop:quantum-limit}
Let $u_j$ be the real normalized eigenfunctions of Proposition~\ref{prop:construction}, and put $h_j=\lambda_j^{-1/2}$. Denote by $\mu_{\sigma\tau}$ the normalized Liouville measure on the graph torus $\Lambda_{\sigma\tau}(0)$ in \eqref{eq:liouville-four-tori}. Then the full sequence has the unique semiclassical measure
\begin{equation}\label{eq:two-component-limit}
 \mu=\tfrac12\mu_{+-}+\tfrac12\mu_{-+}.
\end{equation}
In particular, for every $c\in C_c^\infty(T^*M)$,
\[
 \lim_{j\to\infty}
 \big\langle\operatorname{Op}_{h_j}^w(c)u_j,u_j\big\rangle_{L^2(M)}
 =\frac12\int_{\Lambda_{+-}(0)}c\,d\mu_{+-}
  +\frac12\int_{\Lambda_{-+}(0)}c\,d\mu_{-+}.
\]
The same conclusion holds for any sequence of real normalized functions in $\E_{h_j,0}$, with $h_j\to0$, that vanishes identically on $\gamma_{r,s}$.
\end{proposition}

\begin{proof}
Write $p_L=(\xi_x^2+\xi_y^2)/F$ for the principal symbol of $h^2L$. The joint equations and ellipticity imply that any subsequential defect measure is supported on $\{p_L=1,q=0\}$. Invariance under both commuting Hamiltonian flows gives
\[
 \mu=\sum_{\sigma,\tau=\pm1}c_{\sigma\tau}\mu_{\sigma\tau},
 \qquad c_{\sigma\tau}\ge0,\qquad
 \sum_{\sigma,\tau}c_{\sigma\tau}=1,
\]
as in \eqref{eq:finite-measure}. We show that the common zero curve determines these weights.

Along $\gamma_{r,s}(t)=(st,rt)$ one has
\[
 A_0=r\sqrt{f(rs t)},\qquad B_0=s\sqrt{f(rs t)},
 \qquad \dot\gamma_{r,s}=(s,r).
\]
The four tangential momenta are therefore
\[
 (\sigma A_0,\tau B_0)(\dot\gamma_{r,s})
 =rs\sqrt{f(rs t)}(\sigma+\tau).
\]
The $++$ and $--$ branches each have a unique, nonzero tangential covector, whereas the two mixed branches both restrict to zero. Thus a short arc has tangential windows isolating the $++$ and $--$ sheets from the entire fiber.

In a coordinate neighborhood of this arc, use the defining function $\rho=rx-sy$, up to an additive constant. On the two isolated branches,
\[
 H_{p_L}\rho\big|_{\Lambda_{++}(0)}
 =\frac{2(r^2-s^2)\sqrt{f(rs t)}}{F},
 \qquad
 H_{p_L}\rho\big|_{\Lambda_{--}(0)}
 =-\frac{2(r^2-s^2)\sqrt{f(rs t)}}{F}.
\]
Both are nonzero because $r\ne s$. The first-order symbol $p_1=\sqrt{1+p_L}$ has the same transversality. Proposition~\ref{prop:local-cauchy} therefore supplies cutoffs $\chi_{++}$ and $\chi_{--}$ in flow tubes on the corresponding components such that
\[
 \norm{\operatorname{Op}_{h_j}^w(\chi_{\sigma\sigma})u_j}_{L^2(M)}
 \le C\norm{u_j|_{\gamma_{r,s}}}_{L^2(\gamma_{r,s})}+o(1)
 =o(1),\qquad \sigma=\pm1.
\]
Choose each cutoff nonzero on an open part of its torus and supported away from the other components. Passing to the defect measure gives
\[
 c_{\sigma\sigma}
 \int_{\Lambda_{\sigma\sigma}(0)}|\chi_{\sigma\sigma}|^2
 \,d\mu_{\sigma\sigma}=0.
\]
The integrals are positive, hence $c_{++}=c_{--}=0$.

Since the functions $u_j$ are real, their defect measure is invariant under time reversal $\iota(x,\xi)=(x,-\xi)$. This follows by conjugating a pseudodifferential observable; its principal symbol is changed from $c(x,\xi)$ to $\overline{c(x,-\xi)}$. The map $\iota$ exchanges $\Lambda_{+-}(0)$ and $\Lambda_{-+}(0)$ and sends their normalized Liouville measures to one another. Thus $c_{+-}=c_{-+}$. Their sum is one, proving \eqref{eq:two-component-limit}. Every subsequential defect measure has been identified with the same probability measure, so the entire sequence converges. The proof uses only reality, the joint equations, and the zero trace, which also proves the final assertion.
\end{proof}

\begin{remark}[The spatial limit]
The limit has positive spatial density everywhere, despite exact vanishing on the fixed circle. Set
\[
 Z=\int_{\T^2}\frac{F(x,y)}{A_0(x)B_0(y)}\,dx\,dy.
\]
On each graph torus the absolute Jacobian is
\[
 \left|\det\frac{\partial(p_L,q)}{\partial(\xi_x,\xi_y)}\right|
 =\frac{4A_0B_0}{F}.
\]
Thus the unnormalized joint Liouville density is $F\,dx\,dy/(4A_0B_0)$, and the projection of its normalization is
\[
 \pi_*\mu_{\sigma\tau}
 =\frac1Z\frac{F(x,y)}{A_0(x)B_0(y)}\,dx\,dy.
\]
Consequently,
\begin{equation}\label{eq:constructed-spatial-limit}
 |u_j|^2\,dV_g\ \rightharpoonup\
 \frac1Z\frac{F}{A_0B_0}\,dx\,dy.
\end{equation}
This positive spatial limit is compatible with exact nodality because the two remaining branches have indistinguishable tangential covectors and their traces cancel.
\end{remark}

%% file: sections/counterexample_figures.tex
\subsection{Geometry of the explicit example}

\begin{figure}[H]
 \centering
 \includegraphics[width=.9\textwidth]{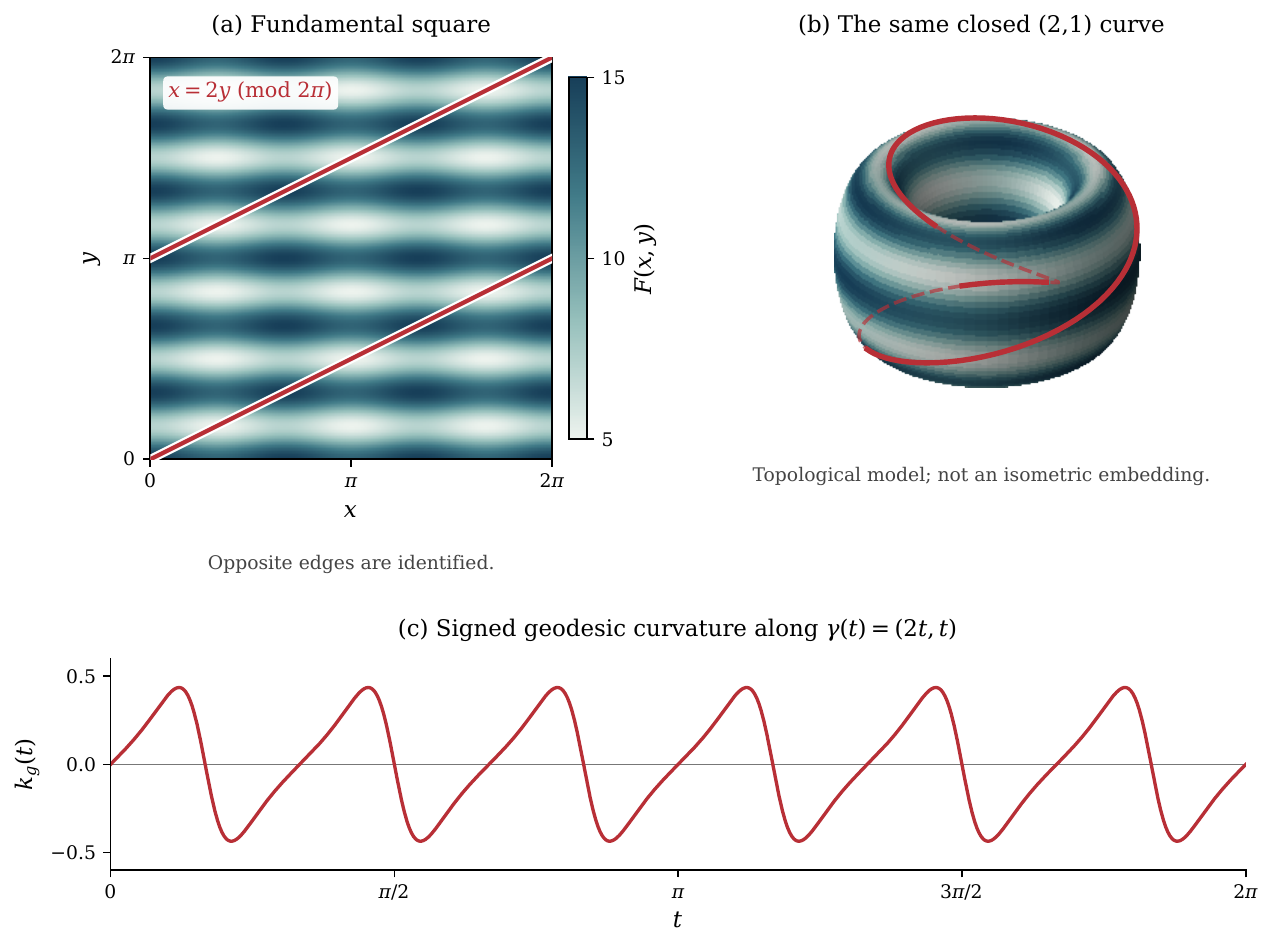}
 \caption{The analytic metric $g=(10+\cos3x+4\cos6y)(dx^2+dy^2)$ and its common nodal circle $\gamma(t)=(2t,t)$. Top left: a coordinate fundamental square; the two red segments join to one circle. Top right: the winding $(2,1)$ on a topological torus. Bottom: the signed geodesic curvature.}
 \label{fig:geometry}
\end{figure}

\begin{figure}[H]
 \centering
 \includegraphics[width=.6\textwidth]{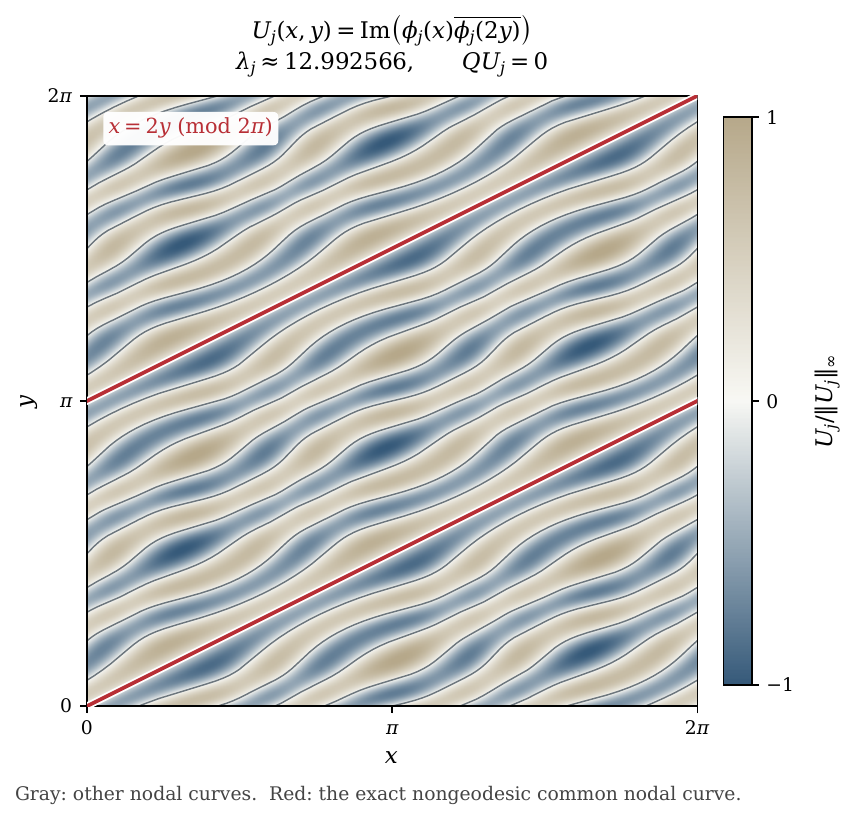}
 \caption{A numerical member of the exact family $U_j(x,y)=\Impart(\phi_j(x)\overline{\phi_j(2y)})$ for $\varepsilon=1$, rescaled for display. The scalar eigenfunction is computed in the twisted Fourier basis $e^{i(3k+1)t}$, with eigenvalue approximately $12.99256648$. }
 \label{fig:eigenfunction}
\end{figure}

For $\varepsilon=1$, the conformal factor ranges from $5$ to $15$. The curvature in \eqref{eq:explicit-curvature} vanishes only at $t=k\pi/6$ and is nonzero on the intervening arcs. Figures~\ref{fig:geometry} and \ref{fig:eigenfunction} display the geometry and one eigenfunction. The three-dimensional drawing records the topology; the metric and curvature are specified by their formulas, not by an embedding in Euclidean space.


%% file: qci_restriction_persistent_nodal_curves.bbl
\begin{thebibliography}{BDWZ12}
\bibitem[BDWZ12]{BDWZ12}
N.~Burq, S.~Dyatlov, R.~Ward, and M.~Zworski,
\emph{Weighted eigenfunction estimates with applications to compressed sensing},
SIAM J. Math. Anal. \textbf{44} (2012), no.~5, 3481--3501.

\bibitem[BGT04]{BGT04}
N.~Burq, P.~G\'erard, and N.~Tzvetkov,
\emph{The Cauchy problem for the nonlinear Schr\"odinger equation on compact manifolds},
in \emph{Phase space analysis of partial differential equations}, Vol.~I,
Pubbl. Cent. Ric. Mat. Ennio Giorgi, Scuola Norm. Sup., Pisa, 2004, 21--52.

\bibitem[BGT07]{BGT07}
N.~Burq, P.~G\'erard, and N.~Tzvetkov,
\emph{Restrictions of the Laplace--Beltrami eigenfunctions to submanifolds},
Duke Math. J. \textbf{138} (2007), no.~3, 445--486.

\bibitem[BR09]{BR09}
J.~Bourgain and Z.~Rudnick,
\emph{Restriction of toral eigenfunctions to hypersurfaces},
C. R. Math. Acad. Sci. Paris \textbf{347} (2009), no.~21--22, 1249--1253.

\bibitem[BR11]{BR11}
J.~Bourgain and Z.~Rudnick,
\emph{On the nodal sets of toral eigenfunctions},
Invent. Math. \textbf{185} (2011), no.~1, 199--237.
\href{https://doi.org/10.1007/s00222-010-0307-2}{doi:10.1007/s00222-010-0307-2};
\href{https://arxiv.org/abs/1003.1743}{arXiv:1003.1743}.

\bibitem[DZ13]{DZ13}
S.~Dyatlov and M.~Zworski,
\emph{Quantum ergodicity for restrictions to hypersurfaces},
Nonlinearity \textbf{26} (2013), no.~1, 35--52.

\bibitem[EPS15]{EPS}
A.~Enciso and D.~Peralta-Salas,
\emph{Eigenfunctions with prescribed nodal sets},
J. Differential Geom. \textbf{101} (2015), no.~2, 197--211.
\href{https://arxiv.org/abs/1404.1039}{arXiv:1404.1039}.

\bibitem[GT20]{GT}
J.~Galkowski and J.~A.~Toth,
\emph{Pointwise bounds for joint eigenfunctions of quantum completely integrable systems},
Comm. Math. Phys. \textbf{375} (2020), 915--947.
\href{https://doi.org/10.1007/s00220-020-03730-3}{doi:10.1007/s00220-020-03730-3}.

\bibitem[GZ21]{GZ21}
J.~Galkowski and S.~Zelditch,
\emph{Lower bounds for Cauchy data on curves in a negatively curved surface},
Israel J. Math. \textbf{244} (2021), no.~2, 971--1000.

\bibitem[Han]{Han}
X.~Han,
\emph{Problems}, online seminar notes, \S4.1.5, Problem~4.9, p.~61, undated.
\url{https://www.csun.edu/~xiaolong/Han_Pblm.pdf}
(accessed 16 October 2026).

\bibitem[JZ16]{JZ16}
J.~Jung and S.~Zelditch,
\emph{Number of nodal domains and singular points of eigenfunctions of negatively curved surfaces with an isometric involution},
J. Differential Geom. \textbf{102} (2016), no.~1, 37--66.

\bibitem[KTZ07]{KTZ07}
H.~Koch, D.~Tataru, and M.~Zworski,
\emph{Semiclassical $L^p$ estimates},
Ann. Henri Poincar\'e \textbf{8} (2007), no.~5, 885--916.

\bibitem[Tot09]{Tot09}
J.~A.~Toth,
\emph{$L^2$-restriction bounds for eigenfunctions along curves in the quantum completely integrable case},
Comm. Math. Phys. \textbf{288} (2009), no.~1, 379--401.

\bibitem[TZ03]{TZ03}
J.~A.~Toth and S.~Zelditch,
\emph{Norms of modes and quasi-modes revisited},
in \emph{Harmonic analysis at Mount Holyoke (South Hadley, MA, 2001)},
Contemp. Math., vol.~320, Amer. Math. Soc., Providence, RI, 2003, 435--458.

\bibitem[TZ13]{TZ13}
J.~A.~Toth and S.~Zelditch,
\emph{Quantum ergodic restriction theorems: manifolds without boundary},
Geom. Funct. Anal. \textbf{23} (2013), no.~2, 715--775.

\bibitem[TZ21]{TZ21}
J.~A.~Toth and S.~Zelditch,
\emph{Nodal intersections and geometric control},
J. Differential Geom. \textbf{117} (2021), no.~2, 345--393.
\href{https://doi.org/10.4310/jdg/1612975018}{doi:10.4310/jdg/1612975018}.
The geodesicity question cited here is in \S1.3, p.~12, of
\href{https://arxiv.org/abs/1708.05754v1}{arXiv:1708.05754v1} (2017).

\bibitem[WW]{WW}
W.~Wang and X.~Wu,
\emph{Average estimate for eigenfunctions along geodesics in the quantum completely integrable case},
\href{https://arxiv.org/abs/2506.15992v2}{arXiv:2506.15992v2} (2026).

\bibitem[WXb]{WXb}
X.~Wu and X.~Xiao,
\emph{Trace estimates and improved pointwise bounds for joint eigenfunctions},
\href{https://arxiv.org/abs/2604.22197}{arXiv:2604.22197} (2026).

\bibitem[Zwo12]{Zwo12}
M.~Zworski,
\emph{Semiclassical analysis},
Graduate Studies in Mathematics, vol.~138, Amer. Math. Soc., Providence, RI, 2012.
\end{thebibliography}
